\documentclass[12pt]{article}
\usepackage[dvips]{color}
\usepackage{graphicx}
\usepackage{epsfig}
\usepackage{verbatim}
\usepackage{amsthm}
\usepackage{amssymb}
\usepackage{amsxtra}     
\usepackage{multirow}

\newtheorem{theorem}{Theorem}[section]
\newtheorem{definition}[theorem]{Definition}
\newtheorem{lemma}[theorem]{Lemma}
\newtheorem{corollary}[theorem]{Corollary}

\newcommand{\PROOF}{\noindent {\bf Proof}: }

\DeclareSymbolFont{AMSb}{U}{msb}{m}{n}
\DeclareMathSymbol{\N}{\mathbin}{AMSb}{"4E}
\DeclareMathSymbol{\Z}{\mathbin}{AMSb}{"5A}
\DeclareMathSymbol{\R}{\mathbin}{AMSb}{"52}
\DeclareMathSymbol{\Q}{\mathbin}{AMSb}{"51}
\DeclareMathSymbol{\I}{\mathbin}{AMSb}{"49}
\DeclareMathSymbol{\C}{\mathbin}{AMSb}{"43}

\begin{document}

\title{The Cyclic Coloring Complex of a $k$-Uniform Hypergraph}
\author{Sarah Crown Rundell\\Denison University} \maketitle

\begin{abstract} In this paper, we study the homology of the cyclic coloring complex of three different types of $k$-uniform hypergraphs.  For the case of a complete $k$-uniform hypergraph, we show that the dimension of the $(n-k-1)^{st}$ homology group is given by a binomial coefficient.  Further, we discuss a complex whose $r$-faces consist of all ordered set partitions $[B_1, ..., B_{r+2}]$ where none of the $B_i$ contain a hyperedge of the complete $k$-uniform hypergraph $H$ and where $1 \in B_1$.  It is shown that the dimensions of the homology groups of this complex are given by binomial coefficients.  As a consequence, this result gives the dimensions of the multilinear parts of the cyclic homology groups of $\C[x_1,...,x_n]/ \{x_{i_1}...x_{i_k} \mid i_{1}...i_{k}$ is a hyperedge of $H \}$. For the other two types of hypergraphs, star hypergraphs and diagonal hypergraphs, we show that the dimensions of the homology groups of their cyclic coloring complexes are given by binomial coefficients as well.
\end{abstract}

\section{Introduction}

\hspace{.5in}In this paper, we will study the homology of the cyclic coloring complex of a $k$-uniform hypergraph.  Let $G$ be a simple graph on $n$ vertices.  The coloring complex of a graph, $\Lambda(G)$ was introduced by Steingr\'{\i}msson~\cite{sm}, and for the case where $G$ has at least one edge, the homology of $\Lambda(G)$ was shown to be nonzero only in dimension $n-3$ by Jonsson~\cite{jo}.  Further, Jonsson showed that in this case the dimension of the nonzero homology group of $\Lambda(G)$ was equal to $\chi_G(-1)-1$, where $\chi_{G}(\lambda)$ is the chromatic polynomial of $G$.  Crown~\cite{cr} defined and studied the cyclic coloring complex of a graph, $\Delta(G)$ and determined the dimensions of the homology groups of $\Delta(G)$.  In particular, she showed that if $G$ has $r$ connected components, the dimension of the $(n-3)^{rd}$ homology group of $\Delta(G)$ is equal to $n-(r+1)$ plus $\frac{1}{r!}\mid \chi_{G}^{r}(0) \mid$, where $\chi_{G}^{r}(\lambda)$ is the $r^{th}$ derivative of $\chi_{G}(\lambda)$.  

\hspace{.5in}Let $H$ be a hypergraph on $n$ vertices.  The coloring complex of a hypergraph, $\Lambda(H)$ was introduced in Long and Rundell~\cite{lr} as well as in Breuer, Dall, and Kubitze~\cite{bd}.  Hanlon~\cite{ha} showed that there exists a Hodge decomposition of the nonzero homology group of $\Lambda(G)$ and that the dimension of the $j^{th}$ Hodge piece of this decomposition equals the absolute value of the coefficient of $\lambda^{j}$ in $\chi_{G}(\lambda)$.  Long and Rundell~\cite{lr} extended Hanlon's result by showing that the Euler Characteristic of the $j^{th}$ Hodge subcomplex of $\Lambda(H)$ is related to the coefficient of $\lambda^{j}$ in $\chi_{H}(\lambda)$.  Long and Rundell also showed that for a class of hypergraphs, which they call star hypergraphs, the coloring complex of the hypergraph is Cohen-Macaulay.  Breuer, Dall, and Kubitzke~\cite{bd} showed that the $f-$ and $h-$ vectors of the coloring complexes of hypergraphs provide tighter bounds on the chromatic polynomials of hypergraphs.  They also show that the coloring complex of a hypergraph has a wedge decomposition, and they provide a characterization of those hypergraphs which have a connected coloring complex.

\hspace{.5in}In this paper, we study the cyclic coloring complex of a complete $k$-uniform hypergraph.  We determine the dimension of the homology group, $HC_{n-k-1}(\Delta(H))$.  We also define a complex $\Delta(H)^C$, whose $r$-faces consist of all ordered set partitions $[B_1, ..., B_{r+2}]$, where none of the $B_i$ contain a hyperedge of $H$ and where $1 \in B_1$.  We compute the dimensions of the homology groups of this complex, $HC_{r}(\Delta(H)^C)$, for $r \geq n-k$.  This gives, as a result, the dimensions of the multilinear parts of the cyclic homology groups of $\C[x_1,...,x_n]/\{x_{i_1}...x_{i_k} \mid i_1...i_k$ is a hyperedge of $H \}$.  Further, for $k = n-1$ and $k=n-2$, we determine the dimensions of $HC_{r}(\Delta(H))$ for all $r$.  In this paper, we will also study the cyclic coloring complex of a star hypergraph and the cyclic coloring complex of a diagonal hypergraph.

\section{Preliminaries}

\begin{definition}
A \emph{hypergraph}, $H$, is an ordered pair, $(V,E)$, where $V$ is a set of vertices and $E$ is a set of subsets of the vertices of $V$.  A hypergraph is said to be \emph{uniform of rank $k$}, or $k$-uniform, if all of its hyperedges have size $k$.  A hypergraph is \emph{complete $k$-uniform} if every subset of size $k$ of $V$ is a hyperedge of $H$.
\end{definition}

Throughout this paper, $H$ will denote a hypergraph whose vertex set $V$ is $\{1,...,n\}$.

\hspace{.5in}Next we define the coloring complex of a hypergraph, $H$, following the presentation in Jonsson~\cite{jo}.

\hspace{.5in}Let $ (B_{1},...,B_{r+2}) $ be an ordered partition
of $ \{ 1,...,n \} $ where at least one of the $B_{i}$ contains a
hyperedge of $H$, and let $\Lambda_{r}$ be the set of ordered partitions
$ ( B_{1},...,B_{r+2} ) $.

\begin{definition} The \emph{coloring complex} of $H$ is the simplicial complex defined by the sequence:

\begin{center}
$...\rightarrow V_{r} \stackrel{\delta_{r}}{\rightarrow} V_{r-1}
\stackrel{\delta_{r-1}}{\rightarrow} ...
\stackrel{\delta_1}{\rightarrow} V_{0}
\stackrel{\delta_0}{\rightarrow} V_{-1}
\stackrel{\delta_{-1}}{\rightarrow} 0$
\end{center}

where $V_{r}$ is the vector space over a field of characteristic
zero with basis $\Lambda_{r}$ and

\begin{center}
$\delta_r( (B_{1},...,B_{r+2}) ) := \displaystyle
\sum_{i=1}^{r+1} (-1)^i ( B_1,...,B_{i} \bigcup B_{i+1},...,
B_{r+2} ) $.
\end{center}

\end{definition}

\hspace{.5in}Notice that $\delta_{r-1} \circ \delta_{r} = 0$.  Then:

\begin{definition}
The $r^{th}$ homology group of $\Lambda(H)$, $H_{r}(\Lambda(H)) :=
ker(\delta_{r})/im(\delta_{r+1})$.
\end{definition}

\hspace{.5in}It is also worth noting that Hultman~\cite{hu} defined a complex that includes both Steingr\'{\i}msson's coloring complex and the coloring complex of a hypergraph as a special case.

\hspace{.5in}To define the cyclic coloring complex, we first must define an equivalence relation on the elements of $\pm \Lambda_r$:

\hspace{.5in}Let $\sigma \in S_{r+2}$ be the $(r+2)$-cycle $(1, 2, ..., r+2)$.  Define $\Delta_r = \pm \Lambda_r / \sim$, where $\sim$ is defined by $(B_1, ..., B_{r+2}) \sim (-1)^{r+1}(B_{\sigma(1)}, ..., B_{\sigma(r+2)})$.  Let $[B_1,...,B_{r+2}]$ denote the equivalence class containing $(B_1, ..., B_{r+2})$.  We will represent each equivalence class of $\Delta_r$ by the unique representative that has $1 \in B_1$.

\hspace{.5in}Let
$$\partial_r([B_1,...,B_{r+2}]) := \sum_{i=1}^{r+1}(-1)^{i+1}[B_1, ..., B_i \cup B_{i+1}, ..., B_{r+2}] + (-1)^{r+3}[B_1 \cup B_{r+2}, B_2,..., B_{r+1}].$$

\hspace{.5in}It is straightforward to check that $\partial$ is well-defined on equivalence classes.  Thus,

\begin{definition}
The \emph{cyclic coloring complex} of $H$, $\Delta(H)$, is the sequence
\begin{center}
$...\rightarrow C_{r} \stackrel{\partial_{r}}{\rightarrow} C_{r-1}
\stackrel{\partial_{r-1}}{\rightarrow} ...
\stackrel{\partial_1}{\rightarrow} C_{0}
\stackrel{\partial_0}{\rightarrow} C_{-1}
\stackrel{\partial_{-1}}{\rightarrow} 0$
\end{center}

where $C_{r}$ is the vector space over a field of characteristic
zero.
\end{definition}

Notice that $\partial_{r-1} \circ \partial_r = 0$, so then:

\begin{definition}
The $r^{th}$ homology group of $\Delta(H)$, $HC_{r}(\Delta(H)) := ker(\partial_r)/im(\partial_{r+1})$.
\end{definition}

\hspace{.5in}As mentioned in Crown~\cite{cr}, the motivation for the definition cyclic coloring complex comes from cyclic homology.  See Loday~\cite{l1} for more information on cyclic homology.  

\hspace{.5in}In a couple of our arguments, we will consider the homology of the quotient of two cyclic coloring complexes, so we will define this quotient now:

\hspace{.5in}Consider the coloring complex of a hypergraph $H$, $\Delta(H)$, and consider a subcomplex, $\Delta(I)$, of $\Delta(H)$.  Then $\Delta(H)/\Delta(I)$ will consist of partitions of $\Delta(H)$ where none of the $B_i$ contains a hyperedge of $I$.  Thus, we obtain the sequence of complexes:
$$\Delta(I) \hookrightarrow \Delta(H) \stackrel{j}{\rightarrow} \Delta(H)/\Delta(I)$$
where $j$ is the quotient map.  From the homology of the pair, $(\Delta(H), \Delta(I))$, this then induces the long exact sequence:
$$\hdots \rightarrow HC_{r}(\Delta(I)) \stackrel{i^{*}}{\rightarrow} HC_{r}(\Delta(H)) \stackrel{j^{*}}{\rightarrow} HC_{r}(\Delta(H)/\Delta(I)) \rightarrow \hdots$$
where $i^{*}$ is the map induced by the inclusion $\Delta(I) \hookrightarrow \Delta(H)$ and $j^{*}$ is the map induced by the quotient map $j$.

\hspace{.5in}One of our results will relate the dimension of $\Lambda(H)$ to the chromatic polynomial of $H$.  So we include the definition of $\chi_{H}(\lambda)$ here:

\begin{definition}
The \emph{chromatic polynomial} of $H$, denoted $\chi_{H}(\lambda)$, is the number of ways to color the vertices of the hyperedges of $H$ with at most $\lambda$ colors, so that the vertices of each edge are colored with at least two colors.
\end{definition}

In this paper, we will also need the following definitions:







\begin{definition} Let $H$ be a uniform hypergraph of rank $k$ with no singleton vertices.  $H$ is a \emph{star hypergraph} if all of its hyperedges intersect in a common set of size $k-1$.
\end{definition}

\begin{definition} Let $H$ be a uniform hypergraph of rank $k$ with edge set $E = \{e_1, ..., e_{n-k+1} \}$.  $H$ is a \emph{diagonal hypergraph} if it isomorphic to the hypergraph with edge set 
$$E = \{ \{1,2,...,k \}, \{2,3,...,k+1\}, \{3,4,...,k+2\}, ..., \{n-k+1, n-k+2,...,n\} \}.$$
\end{definition}

\section{The Complete $k$-Uniform Hypergraph}

\hspace{.5in} In this section, we discuss the homology of the cyclic coloring complex of a complete $k$-uniform hypergraph.  We let $\Delta(E_n)$ be the cyclic coloring complex of the complete graph with looped edges at each vertex.  Notice that then $\Delta_{n-2}(E_n)$ consists of all ordered partitions $[B_1,...,B_{n}]$ of $[n]$ where $\mid B_i \mid = 1$ for all $i$, $1 \leq i \leq n$.  

\begin{theorem}\label{Vertex1Hypergraph}Let $H$ be the $k$-uniform hypergraph on $n$ vertices with edge set consisting of all possible hyperedges of size $k$ containing the vertex 1.  Then the dimension of $HC_r(\Delta(H))$ is nonzero for $n-k-1 \geq r \geq -1$ and is given by
$$\dim HC_{r}(\Delta(H)) = \binom{n-1}{r+1}.$$
\end{theorem}

\PROOF

Consider the function 
$$f([B_1, ..., B_{l+2}]) = \mid B_1 \mid$$
where $B_1$ is the block containing 1.  Let $\Delta_{l}^{m}(H)$ denote the chains in $\Delta_l(H)$ where 
$$f([B_1, ..., B_{l+2}]) = m,$$
and let $\Delta^{m}(H)$ denote the complex formed by the chains in $\Delta_{i}^{m}(H)$, $-1 \leq i \leq n-2$.  Notice that $f$ gives a filtration of each $\Delta_{l}(H)$ and that the boundary map $\partial$ respects the filtration.  We will use a spectral sequence argument to determine the dimensions of the homology groups of $\Delta(H)$.  See Chow~\cite{ch} for an introduction to spectral sequences, and we will follow his notation throughout the proof.

\hspace{0.5in}Notice that $\Delta^{m}(E_n) = \Delta^{m}(H)$ for $n-k+1 \leq m \leq n$.  By the same argument used in the proof of Theorem 3.2 in Crown~\cite{cr}, we can see that the chains of $\Delta^{m}(H)$ can be partitioned into subcomplexes determined by the elements of $B_1$.  The number of subcomplexes of $\Delta^{m}(H)$ is given by the number of ways of forming a subset of size $m-1$ from a set of size $n-1$.  In particular, there are $\binom{n-1}{m-1}$ subcomplexes of $\Delta^{m}(H)$.  Further, as argued in the proof of Theorem 3.2 in Crown ~\cite{cr}, the spectral sequence collapses, and
$$\dim(HC_{r}(\Delta(H)) = \binom{n-1}{r+1}.$$
\qed

\hspace{.5in}Further, the proof of Theorem 3.2 in Crown ~\cite{cr} also shows that the homology representatives of the $r^{th}$ homology group of $\Delta(H)$ are indexed by the subsets of size $n-r-2$ of $\{ 2, ..., n \}$.  Namely, for each subset, $A$, of size $n-r-2$ of $\{2,...,n\}$, we obtain one homology representative of $HC_{r}(\Delta(H))$, $\displaystyle \sum_{\sigma \in S_{r+1}} sgn(\sigma) [ A \cup {1}, a_{\sigma(1)}, ..., a_{\sigma(r+1)}]$, where  $\{ a_1, ..., a_{r+1} \}$ is the complement of $A$ in $\{2, ..., n \}$.

\hspace{.5in}Let $H$ be the complete $k$-uniform hypergraph on $n$ vertices.  Notice that $\Delta(H)$ is a subcomplex of $\Delta(E_n)$.  Let $\Delta(H)^C = \Delta(E_n)/\Delta(H)$.  Also, notice that the boundary map of $\Delta(H)^C$ maps partitions $[B_1,...,B_{r+2}]$ where at least one of the $B_i$ contains an edge of $H$ of zero.  To compute the homology of $\Delta(H)$, we will first compute the homology of $\Delta(H)^C$.

\begin{theorem}\label{Complement}Let $H$ be the complete $k$-uniform hypergraph on $n$ vertices.  

For $n-2 \geq r > n-k$,
$$\dim HC_{r}(\Delta(H)^C) = \binom{n-1}{r+1} $$
and
$$\dim HC_{n-k}(\Delta(H)^C) = \binom{n-1}{n-k-1} + \binom{n-1}{n-k+1}.$$
\end{theorem}

\PROOF

For the first part of the theorem, consider the following long exact sequence:

\begin{center}
\begin{tabular}{ccccccccc}
$0$ & $\rightarrow$ & $HC_{n-2}(\Delta(H))$ &
$\rightarrow$ & $HC_{n-2}(\Delta(E_{n}))$ &
$\rightarrow$ &
$HC_{n-2}(\Delta(H)^C)$ & $\rightarrow$ & \\
& & $HC_{n-3}(\Delta(H))$ & $\rightarrow$ &
$HC_{n-3}(\Delta(E_n)$ &
$\rightarrow$ & $HC_{n-3}(\Delta(H)^C)$ & $\rightarrow$ &  \\
&  & $HC_{n-4}(\Delta(H))$ & $\rightarrow$ &
$HC_{n-4}(\Delta(E_n))$ & $\rightarrow$ &
$HC_{n-4}(\Delta(H)^C)$ & $\rightarrow$ &
\end{tabular}
\newline \centerline{$\vdots$} \newline
\end{center}

Notice that since $HC_{r}(\Delta(H)) = 0$ for $r > n-k-1$, by exactness, $\dim HC_{r}(\Delta(H)^C) = \dim HC_{r}(\Delta(E_n)) = \binom{n-1}{r+1}$ for $n-2 \geq r > n-k$.

\hspace{.5in}We will now prove the second statement of the theorem.  Let $W_1, ..., W_{\binom{n}{k}}$ be the subsets of $\{1, ..., n \}$ of size $k$, listed in lexicographic order.  Let $\Delta(E_{n}^{(0, W_l)})$ be the complex formed be the chains $[B_1, ..., B_{r+2}]$ where for all $i$, $1 \leq i \leq l$, the elements of $W_i$ are not in the same block of the partition.  Let $\Delta(E_{n}^{(1, W_l)})$ be the complex formed by the chains $[B_1, ..., B_{r+2}]$ where for all $i$, $1\leq i \leq l-1$, the elements of $W_i$ are not in the same block of the partition, but the elements of $W_l$ are in the same $B_j$, for some $j$, $1 \leq j \leq r+2$.  Notice then that:
$$\Delta(E_n^{(0,W_{l-1})})/\Delta(E_n^{(1,W_l)}) = \Delta(E_n^{(0,W_l)}) .$$
Using this notation, $\Delta(H)^C = \Delta(E_n^{(0,W_{\binom{n}{k}})})$.  We will compute the homology of $\Delta(H)^C$ by sequentially computing the homology of $\Delta(E_n^{(0,W_l)})$.

\hspace{.5in}  Let $H'$ be the $k$-uniform hypergraph on $n$ vertices with edge set consisting of all possible hyperedges of size $k$ containing the vertex 1.  Since $1 \in W_l$ for $1 \leq l \leq \binom{n-1}{k-1}$, $\Delta(H')^C = \Delta(E_n^{(0,W_{\binom{n-1}{k-1}})})$.  We begin by using the following long exact sequence to calculate the homology of $\Delta(H')^C$:

\begin{center}
\begin{tabular}{ccccccccc}
$0$ & $\rightarrow$ & $HC_{n-2}(\Delta(H'))$ &
$\stackrel{\alpha_{n-2}}{\rightarrow}$ & $HC_{n-2}(\Delta(E_{n}))$ &
$\rightarrow$ &
$HC_{n-2}(\Delta(H')^C)$ & $\rightarrow$ & \\
& & $HC_{n-3}(\Delta(H'))$ & $\stackrel{\alpha_{n-3}}{\rightarrow}$ &
$HC_{n-3}(\Delta(E_n)$ &
$\rightarrow$ & $HC_{n-3}(\Delta(H')^C)$ & $\rightarrow$ &  \\
&  & $HC_{n-4}(\Delta(H'))$ & $\stackrel{\alpha_{n-4}}{\rightarrow}$ &
$HC_{n-4}(\Delta(E_n))$ & $\rightarrow$ &
$HC_{n-4}(\Delta(H')^C)$ & $\rightarrow$ &
\end{tabular}
\newline \centerline{$\vdots$} \newline
\end{center}

By the same argument as in the first part of this proof, we can see that the dimension of $HC_{r}(\Delta(H'))^C$ for $n-2 \geq r \geq n-k+1$ is $\binom{n-1}{r+1}$.  So consider $r = n-k$.  We noted after the proof of Theorem~\ref{Vertex1Hypergraph} that the homology representatives of the $r^{th}$ homology group of $\Delta(H')$ are indexed by the subsets of size $n-r-2$ of $\{ 2, ..., n \}$.  As noted in the proof of Theorem 3.2 in Crown~\cite{cr}, these are the same as the homology representatives of the $r^{th}$ homology group of $\Delta(E_n)$.  Thus $\alpha_r$ is bijective for all $r \leq n-k-1$.   By the above argument, it suffices to consider the exact sequence:

\begin{center}
\begin{tabular}{ccccccc}
$0$ & $\rightarrow$ & $HC_{n-k}(\Delta(E_n))$ & $\rightarrow$ & 
$HC_{n-k}(\Delta(H')^C)$ & $\stackrel{\phi_{n-k}}{\rightarrow}$ & \\ 
$HC_{n-k-1}(\Delta(H'))$ & $\stackrel{\alpha_{n-k-1}}{\rightarrow}$ & $HC_{n-k-1}(\Delta(E_n))$ & $\rightarrow$ & 0 & & \\
\end{tabular}
\end{center}

Since $\alpha_{n-k-1}$ is injective, the image of $\phi_{n-k}$ is zero.  Therefore, the dimension of $HC_{n-k}(\Delta(H')^C)$ equals the dimension of the kernel of $\phi_{n-k}$.  By exactness, the dimension of the kernel of $\phi_{n-k}$ equals the dimension of $HC_{n-k}(\Delta(E_n))$.  So the dimension of $HC_{n-k}(\Delta(H')^C)$ is $\binom{n-1}{n-k+1}$.



\hspace{.5in}Now we will sequentially compute $HC_{r}(\Delta(H)^C)$.  Let $l > \binom{n-1}{k-1}$.  Notice that $HC_{r}(\Delta(E_n^{(1,W_l)})$ is zero for $r \geq n-k$.  So by exactness, $\dim HC_{r}(\Delta(E_n^{(0,W_{l})}) = \dim HC_{r}(\Delta(H')^C) = \binom{n-1}{r+1}$ for $r > n-k$.  Further, we know $HC_{r}(\Delta(H')^C) = 0$ for $r \leq n-k-1$.  It suffices then to consider the the exact sequence:

\begin{center}
\begin{tabular}{ccccccc}
0 & $\rightarrow$ & $HC_{n-k}(\Delta(H')^C)$ & $\rightarrow$ & $HC_{n-k}(\Delta(E_n^{(0,W_{\binom{n-1}{k-1}+1})}))$ & $\rightarrow$ \\
 $HC_{n-k-1}(\Delta(E_n^{(1,W_{\binom{n-1}{k-1}+1})}))$ &
$\stackrel{\alpha_{n-k-1}}{\rightarrow}$ & 0 &
$\rightarrow$ &
$HC_{n-k-1}(\Delta(E_n^{(0,W_{\binom{n-1}{k-1}+1})}))$ & $\rightarrow$ & \\
$HC_{n-k-2}(\Delta(E_n^{(1,W_{\binom{n-1}{k-1}+1})}))$ &
$\stackrel{\alpha_{n-k-2}}{\rightarrow}$ & 0 &
$\rightarrow$ &
$HC_{n-k-2}(\Delta(E_n^{(0,W_{\binom{n-1}{k-1}+1})}))$ & $\rightarrow$ & \\
\end{tabular}
\newline \centerline{$\vdots$} \newline
\end{center}

By exactness, $HC_{r}(\Delta(E_n^{(0,W_{\binom{n-1}{k-1}+1})})) \cong HC_{r-1}(\Delta(E_n^{(1,W_{\binom{n-1}{k-1}+1})}))$ for $ n-k-1 \geq r \geq 0$.  Consider the complex $\Delta(E_n^{(1,W_{\binom{n-1}{k-1}+1})})$.  Notice that the homology of this complex is equal to the homology of the complex $\Delta(T_{n-k+1}^{(0,1))})$, where $T_{n-k+1}$ is the tree on $n-k+1$ vertices, with root labeled 1, and edges $(1,2)$, $(1,3)$, ..., $(1, n-k+1)$, and where vertices $1$ and $2$ are not in the same block of a chain $[B_1, ..., B_{r+2}]$.  Let $\{ a_1, ..., a_{n-k-1} \}$ be the elements of the complement of $W_{\binom{n-1}{k-1}+1}$ in $\{ 2, ..., n \}$ listed in increasing order.  The isomorphism is given by mapping 1 to 1, $W_{\binom{n-1}{k-1}+1}$ to 2, $a_1$ to 3, ...., $a_{n-k-1}$ to $n-k+1$.  By Lemma 3.3 of Crown~\cite{cr}, the dimension of $HC_{n-k-1}(\Delta(E_n^{(1,W_{\binom{n-1}{k-1}+1})}))$ is equal to $\binom{(n-k+1)-(1+1)}{((n-k-1)+2)-(1+1)} = 1$, and therefore, $\dim HC_{n-k}(\Delta(E_n^{(0,W_{\binom{n-1}{k-1}+1})})) = 1 + \binom{n-1}{n-k+1}$.
 
 \hspace{.5in}We claim that for each subsequent set $W_l$ removed, the contribution to $\dim HC_{n-k}(\Delta(H)^C)$ will be one.  Notice that for all $l$, $ \binom{n-1}{k-1} + 1 \leq l \leq \binom{n}{k}$, the homology of the complex $\Delta(E_n^{(1, W_l)})$ is equal to the homology of the complex $\Delta(T_{n-k+1}^{(0,i)})$ for some $i$.  By Lemma 3.3 of Crown, the dimension of $HC_{n-k-1}(\Delta(E_n^{(0,W_l)}))$ is equal to $\binom{(n-k+1)-(l+1)}{((n-k-1)+2)-(l+1)} = 1$.  By exactness, for each subsequent set $W_l$ removed, the contribution to $\dim HC_{n-k}(\Delta(H)^C)$ will be one.  Since there are $\binom{n-1}{k}$ sets $W_l$ of size $k$ of $[n-1]$, $\dim HC_{n-k}(\Delta(H)^C) = \binom{n-1}{k} + \binom{n-1}{n-k+1} = \binom{n-1}{n-k-1} + \binom{n-1}{n-k+1}$. 
\qed

\hspace{.5in}Notice that the elements of $\Delta_r(H)^C$ are in bijection with the cyclic words $[D_1,...,D_{r+2}]$ where $D_i \in \C[x_1,...,x_n]/\{ x_{i_1}...x_{i_k} \mid i_1...i_k $ is a hyperedge of $H \}$ and $[D_1,...,D_{r+2}]$ is an ordered partition of $x_1...x_n$ with $x_1 \in D_1$.  It then follows that we have the following corollary:

\begin{corollary} For the complete $k$-uniform hypergraph on $n$ vertices, $H$, the dimension of the multilinear part of the $r^{th}$ cyclic homology group of $\C[x_1,...,x_n]/\{ x_{i_1}...x_{i_k} \mid i_1...i_k $ is a hyperedge of $H \}$ is $\binom{n-1}{r+1}$ for $n-k \leq r \leq n-2$ and $\binom{n-1}{n-k-1} + \binom{n-1}{n-k+1}$ for $r = n-k$.
\end{corollary}

We will now determine the dimension the $(n-k-1)^{st}$ homology group of $\Delta(H)$ for a complete $k$-uniform hypergraph:

\begin{theorem}\label{MainResult}
Let $H$ be a complete $k$-uniform hypergraph.  Then
$$\dim HC_{n-k-1}(\Delta(H)) = \binom{n}{n-k}.$$
\end{theorem}

\PROOF

Consider the following long exact sequence:

\begin{center}
\begin{tabular}{cccccccc}
0 & $\rightarrow$ & $HC_{n-k}(\Delta(E_n))$ & $\rightarrow$ & $HC_{n-k}(\Delta(H)^C)$ & $\stackrel{\phi_{n-k}}{\rightarrow}$ & \\
$HC_{n-k-1}(\Delta(H))$ & $\stackrel{\alpha_{n-k-1}}{\rightarrow}$ & $HC_{n-k-1}(\Delta(E_n))$ & $\stackrel{\beta_{n-k-1}}{\rightarrow}$ & $HC_{n-k-1}(\Delta(H)^C)$ & $\rightarrow$ & $\hdots$ \\
\end{tabular}
\end{center}

By the proof of Theorem 3.2 in Crown~\cite{cr}, each of the homology representatives of $HC_{r}(\Delta(E_n))$ corresponds to a subset, $A$, of $\{2, ..., n \}$ of size $n-r-2$.  Since the set $\{1\} \cup A$ is an edge of $H$, each of the homology representatives of $HC_{n-k-1}(\Delta(E_n))$ is mapped to zero by the map $\beta_{n-k-1}$.  Therefore, the dimension of the kernel of $\beta_{n-k-1}$ is $\binom{n-1}{n-k}$.  By exactness and Theorem 3.2 in Crown~\cite{cr}, the dimension of the kernel of $\phi_{n-k} = \binom{n-1}{n-k+1}$.  Thus,
\begin{center}
\begin{eqnarray*}
\dim HC_{n-k-1}(\Delta(H)) & = & \left( \dim HC_{n-k}(\Delta(H)^C) - \binom{n-1}{n-k+1} \right) + \binom{n-1}{n-k} \\
 & = & \binom{n-1}{n-k-1} + \binom{n-1}{n-k} \\
 & = & \binom{n}{n-k} \\
 \end{eqnarray*}
\end{center}

\qed

When $k = n-1$ and $k=n-2$, we have the results:

\begin{theorem}Let $H$ be the complete $(n-1)$-uniform hypergraph on $n$ vertices.  Then
$$\dim HC_{0}(\Delta(H)) = \binom{n}{1} = n$$
and
$$\dim HC_{-1}(\Delta(H)) = \binom{n}{0} = 1.$$
\end{theorem}

\PROOF

Note that there $\binom{n}{1}$ elements in $\Delta_0(H)$ and each of these elements is mapped to zero under $\partial_0$.  Thus, $\dim HC_0(\Delta(H)) = n$.  Since $\Delta_{-1}(H) = \{ [12 ... n ] \}$, it follows that $\dim HC_{-1}(\Delta(H)) = 1$.

\qed



\begin{theorem}Let $H$ be the complete $(n-2)$-uniform hypergraph on $n$ vertices.  Then for $ -1 \leq r \leq 1$, 
$$\dim HC_{r}(\Delta(H)) = \binom{n}{r+1}.$$
\end{theorem}

\PROOF

From Theorem ~\ref{MainResult}, we know that $\dim HC_{1}(\Delta(H)) = \binom{n}{2}$.   The set $\Delta_{1}(H)$ consists of all ordered partitions $ [B_1, B_2, B_3] $ where $\{1\} \in B_1$ and where one of the $B_i$ is a hyperedge of $H$.  It follows then that there are $2 \binom{n-1}{n-3} + 2 \binom{n-1}{n-2} = 2 \binom{n}{n-2}$ elements in $\Delta_{1}(H)$.  Thus the dimension of the image of $\partial_1$ is $\binom{n}{2}$.  The dimension of the kernel of $\partial_0$ equals the cardinality of $\Delta_0(H)$.  Since $\Delta_{0}(H)$ consists of all ordered partitions $[B_1, B_2]$ where $\{1\} \in B_1$ and where one of the $B_i$ contains a hyperedge of $H$, there are $ n + \binom{n}{2}$ elements in $\Delta_0(H)$.  So, the dimension of $HC_{0}(\Delta(H)) = \binom{n}{1}$.  It is clear that the dimension of $HC_{-1}(\Delta(H)) = \binom{n}{0}$.




\qed



\section{Star Hypergraphs}

\hspace{0.5in}  In this section, we will determine the dimensions of the homology groups of the coloring complexes of the star hypergraphs. 

\begin{theorem}  Let $H$ be a star hypergraph on $n$ vertices.  For $n-k-1 \geq r \geq -1$,
$$\dim HC_{r}(\Delta(H)) = \binom{n-k+1}{r+1}.$$
\end{theorem}

\PROOF

\hspace{0.5in}Without loss of generality, label the vertices of the common intersection set $\{1, ..., k-1 \}$.  Notice that since the hyperedges of $H$ intersect in a common set of size $k-1$, the cyclic coloring complex of $H$ is equivalent to the cyclic coloring complex of a star graph with vertices $\{ 12\hdots k-1, k, ..., n \}$ and edges $\{12\hdots k-1,k\}, \{12\hdots k-1,k+1\}, ..., \{12\hdots k-1,n\}$.  By Corollary 3.4 in Crown~\cite{cr}, the dimension of the $r^{th}$ homology group of the cyclic coloring complex of this star graph is $\displaystyle \binom{(n-k+2)-1}{r+1}$.
\qed

\begin{corollary}  Let $H$ be a hypergraph on $n$ vertices consisting of a star hypergraph (with at least two edges of size $k$) on $l$ of the vertices, together with $n-l$ singleton vertices.  Let $G$ be the graph on $n-k+1$ vertices obtained by contracting the common set of size $k-1$ of the edges of the star hypergraph.  Then 
$$\dim HC_{n-k-1}(\Delta(H)) = \frac{1}{(n-l+1)!} \mid \chi_{G}^{(n-l+1)}(0) \mid + (l-k-1).$$
For $r \geq l-k-2$, the dimension of the $r^{th}$ homology group of $\Delta(H)$ is 
$$\binom{n-k}{r+1}-\binom{n-l+1}{(n-k-1)-r}+\binom{n-l}{(n-k-1)-(r+1)}\frac{1}{(n-l+1)!}\mid \chi_{G}^{(n-l+1)}(0) \mid.$$
and for $r < l-k-2$, the dimension of the $r^{th}$ homology group of $\Delta(H)$ is
$$\binom{n-k}{r+1}.$$
\end{corollary}

\PROOF
The result follows by noting that $\Delta(H)$ is equivalent to the cyclic coloring complex of $G$ and using Theorem 5.6 of Crown~\cite{cr}.
\qed

\section{Diagonal Hypergraphs}

\hspace{0.5in}In this section, we will determine the dimensions of the homology groups of $\Delta(H)$, where $H$ is a diagonal hypergraph on $n$ vertices.  As in the case of a complete $k$-uniform hypergraph, we will need to first compute the dimensions of the homology groups of $\Delta(H)^C$.

\begin{theorem} Let $H$ be a diagonal hypergraph with $k > \lceil{\frac{n}{2}}\rceil$.  For $r > n-k-1$, 
$$\dim HC_{r}(\Delta(H)^C) = \binom{n-1}{r+1},$$
and for $n-k-1 \geq r \geq 0$,
$$\dim HC_{r}(\Delta(H)^C) = \binom{n-1}{r+1} - (n-k)\binom{n-k-1}{r} - \binom{n-k}{r+1}$$
where $l$ is the number of hyperedges of $H$.
\end{theorem}

\hspace{0.5in}Let $\Delta(E_n^{(1,1)})$ be the subcomplex of $\Delta(E_n)$ formed by the partitions $[B_1,...,B_{r+2}]$ of $\Delta_r(E_n)$ where vertices 2 and 3 are elements of $B_i$ for some $i$.  We then define $\Delta(E_n^{(0,1)})$ to be the complex $\Delta(E_n)/\Delta(E_n^{(1,1)})$, and it is formed by the partitions $[B_1,...,B_{r+2}]$ of $\Delta_r(E_n)$ where vertices 2 and 3 are not elements of the same $B_i$ for some $i$.  To prove the above theorem, we will need the following lemma:

\begin{lemma}The dimension of the $r^{th}$ homology group of $\Delta(E_n^{(0,1)})$ is $\binom{n-2}{r}$.
\end{lemma}

\PROOF

\hspace{0.5in}The proof will follow the base case of Lemma 3.3 of Crown~\cite{cr}.  Consider the following long exact sequence:
\begin{center}
\begin{tabular}{ccccccc}
0 & $\rightarrow$ & $HC_{n-2}(\Delta(E_n))$ & $\rightarrow$ & $HC_{n-2}(\Delta(E_n^{(0,1)}))$ & $\rightarrow$ \\
 $HC_{n-3}(\Delta(E_n^{(1,1)}))$ &
$\stackrel{\alpha_{n-3}}{\rightarrow}$ & $HC_{n-3}(\Delta(E_n))$ &
$\rightarrow$ &
$HC_{n-3}(\Delta(E_n^{(0,1)}))$ & $\rightarrow$ & \\
$HC_{n-4}(\Delta(E_n^{(1,1)}))$ &
$\stackrel{\alpha_{n-4}}{\rightarrow}$ & $HC_{n-4}(\Delta(E_n))$ &
$\rightarrow$ &
$HC_{n-4}(\Delta(E_n^{(0,1)}))$ & $\rightarrow$ & \\
& & $\vdots$ & & & \\
\end{tabular}
\end{center}
Under a relabeling, notice that for all $r$, the elements of $\Delta_r(E_n^{(1,1)})$ are in bijection with the elements of $\Delta_r(E_{n-1})$.  Namely, for a partition $[B_1,...,B_{r+2}]$ in $\Delta_r(E_n^{(1,1)})$ map 1 to 1, 23 to 2, 4 to 3, ..., $n$ to $n-1$.  Thus, the dimension of $HC_{r}(\Delta(E_n^{(1,1)}))$ is $\binom{n-2}{r+1}$.  Further, under the same relabeling, we obtain a set of homology representatives for $HC_{r}(\Delta(E_n^{(1,1)}))$ from the set of homology representatives for $HC_{r}(\Delta(E_{n-1}))$.  Since the homology representatives of $HC_{r}(\Delta(E_{n-1}))$ are indexed by the subsets of size $(n-2) - (r+1)$ of the set $\{2, ..., n-1\}$, the homology representatives of $HC_r(\Delta(E_n^{(1,1)}))$ are indexed by the subsets of size $(n-2) - (r+1)$ of the set $\{23, 4, ...,n\}$.  Since these later representatives are a subset of the homology representatives for $HC_r(\Delta(E_n))$, the map $\alpha_r$ is injective for all $r$.  By exactness, it then follows that
the dimension of $HC_{r}(\Delta(E_n^{(0,1)})$ is $\binom{n-2}{r}$.

\qed

\hspace{0.5in}It is worth noting that the proof of the lemma gives a set of homology representatives of $HC_{r}(\Delta(E_n^{(0,1)}))$.  Namely a set of homology representatives of $HC_{r}(\Delta(E_n^{(0,1)}))$ is given by the subset of homology representatives of $HC_{r}(\Delta(E_n))$ where 2 and 3 are not in the same block of a partition.  It is also worth noting that if we generalize the definition of $\Delta(E_n^{(1,1)})$ to be the subcomplex of $\Delta(E_n)$ formed by the partitions $[B_1,...,B_{r+2}]$ of $\Delta(E_n)$ where vertices $i-1$ and $i$ are in the same block of a partition (for some $i$, $2 \leq i \leq n$), then the above lemma is true.

We now present the proof of Theorem 5.1:

\PROOF

\hspace{0.5in}We will use a similar proof technique to that used in the proof of Theorem~\ref{Complement} above.  Without loss of generality, we will suppose that the edge set of $H$ is 
$$E = \{ \{1,2,...,k \}, \{2,3,...,k+1\}, \{3,4,...,k+2\}, ..., \{n-k+1, n-k+2,...,n\} \}$$ 
and let $e_i$ be the edge $\{ i, i+1, ..., k+(i-1) \}$.  We let $\Delta(E_n^{(1,e_i)})$ be the subcomplex of $\Delta(E_n)$ consisting of all partitions $[B_1, ..., B_{r+2}]$ of $[n]$ such that at least one of the $B_j$ contains the edge $e_i$, $1 \in B_1$, and none of the $B_j$ contain the edges $e_1, ..., e_{i-1}$.  We then define $\Delta(E_n^{(0,e_i)})$ to be the complex $\Delta(E_n)/\Delta(E_n^{(1,e_i)})$.  Notice then that $\Delta(H)^C$ is equal to the complex $\Delta(E_n^{(0,e_{n-k+1})})$.  To determine the dimensions of the homology groups of $\Delta(H)^C$, we will compute the dimensions of the homology groups of $\Delta(E_n^{(0,e_i)})$ for each $i$, $1 \leq i \leq n-k+1$.

\hspace{0.5in}We first compute the dimensions of the homology groups of $\Delta(E_n^{(0,e_1)})$.  Consider the exact sequence:
\begin{center}
\begin{tabular}{ccccccc}
0 & $\rightarrow$ & $HC_{n-2}(\Delta(E_n))$ & $\rightarrow$ & $HC_{n-2}(\Delta(E_n^{(0,e_1)}))$ & $\rightarrow$ \\
& & $\vdots$ & & & \\
0 & $\rightarrow$ & $HC_{n-k}(\Delta(E_n))$ & $\rightarrow$ & $HC_{n-k}(\Delta(E_n^{(0,e_1)}))$ & $\rightarrow$ \\
 $HC_{n-k-1}(\Delta(E_n^{(1,e_1)}))$ &
$\stackrel{\alpha_{n-k-1}}{\rightarrow}$ & $HC_{n-k-1}(\Delta(E_n))$ &
$\rightarrow$ &
$HC_{n-k-1}(\Delta(E_n^{(0,e_1)}))$ & $\rightarrow$ & \\
$HC_{n-k-2}(\Delta(E_n^{(1,e_1)}))$ &
$\stackrel{\alpha_{n-k-2}}{\rightarrow}$ & $HC_{n-k-2}(\Delta(E_n))$ &
$\rightarrow$ &
$HC_{n-k-2}(\Delta(E_n^{(0,e_1)}))$ & $\rightarrow$ & \\
& & $\vdots$ & & & \\
\end{tabular}
\end{center}
Notice that for $n-2 \geq r > n-k$, $\Delta(HC_{r}(\Delta(E_n)) \cong \Delta(HC_{r}(\Delta(E_n^{(0,e_1)}))$.  Consider the complex $\Delta(E_n^{(1,e_1)})$.  If we relabel the vertex 1 in $E_{n-k+1}$ with the elements of $e_1$, then the complex $\Delta(E_n^{(1,e_1)})$ is equivalent to the complex $\Delta(E_{n-k+1})$.  Therefore, for $n-k-1 \geq r \geq -1$,
$$\dim HC_{r}(\Delta(E_{n-k+1}^{(1,e_1)})) = \binom{n-k}{r+1},$$
and further, from the proof of Theorem 3.2 in Crown~\cite{cr}, for each subset, $A$, of size $n-k-1-r$ of $\{ k+1, ..., n-k+1 \}$, we obtain one homology representative of $HC_{r}(\Delta(E_{n-k+1}^{(1,e_1)}))$, $\displaystyle \sum_{\sigma \in S_{r+1}} sgn(\sigma)[ e_1 \cup A, a_{\sigma(1)}, ..., a_{\sigma(r+1)}]$, where $\{ a_1, ..., a_{r+1} \}$ is the complement of $A$ in $\{ k+1, ..., n-k+1 \}$.  Since these are a subset of the homology representatives of $HC_{r}(\Delta(E_n))$, this then implies that the map $\alpha_r$ is injective for all $r$, $n-k-1 \geq r \geq -1$.  By exactness, 
$$\dim HC_{n-k}(\Delta(E_n^{(0,e_1)})) = \binom{n-1}{n-k+1}$$
and
$$\dim HC_{r}(\Delta(E_n^{(0,e_1)})) = \binom{n-1}{r+1} - \binom{n-k}{r+1}.$$ 
Notice that the set of homology representatives of $HC(\Delta(E_n^{(0,e_1)}))$ is the set of homology representatives of $HC(\Delta(E_n))$ where the elements of $e_1$ are not all in the same block of a partition.

\hspace{0.5in}Now consider the long exact sequence
\begin{center}
\begin{tabular}{ccccccc}
0 & $\rightarrow$ & $HC_{n-2}(\Delta(E_n^{(0,e_1)}))$ & $\rightarrow$ & $HC_{n-2}(\Delta(E_n^{(0,e_2)}))$ & $\rightarrow$ \\
& & $\vdots$ & & & \\
0 & $\rightarrow$ & $HC_{n-k}(\Delta(E_n^{(0,e_1)}))$ & $\rightarrow$ & $HC_{n-k}(\Delta(E_n^{(0,e_2)}))$ & $\rightarrow$ \\
 $HC_{n-k-1}(\Delta(E_n^{(1,e_2)}))$ &
$\stackrel{\alpha_{n-k-1}}{\rightarrow}$ & $HC_{n-k-1}(\Delta(E_n^{(0,e_1)}))$ &
$\rightarrow$ &
$HC_{n-k-1}(\Delta(E_n^{(0,e_2)}))$ & $\rightarrow$ & \\
$HC_{n-k-2}(\Delta(E_n^{(1,e_2)}))$ &
$\stackrel{\alpha_{n-k-2}}{\rightarrow}$ & $HC_{n-k-2}(\Delta(E_n^{(0,e_1)}))$ &
$\rightarrow$ &
$HC_{n-k-2}(\Delta(E_n^{(0,e_2)}))$ & $\rightarrow$ & \\
& & $\vdots$ & & & \\
\end{tabular}
\end{center}

By exactness, it follows that for $n-2 \geq r > n-k$, $\dim HC_{r}(\Delta(E_n^{(0,e_2)})) = \binom{n-1}{r+1}$.  Consider the complex $\Delta(E_n^{(1,e_2)})$ and let $T_{n-k+1}$ denote the tree on $n-k+1$ vertices with edges $\{1, 2\}, \{1, 3\}, ..., \{1, n-k+1\}$.  If we relabel vertex 2 in $T_{n-k+1}$ with the elements of $e_2$, vertex 3 with $k+2$, ... , and vertex $n-k+1$ with $n$, then the complex $\Delta(E_n^{(1,e_1)})$ is equivalent to the complex $\Delta(T_{n-k+1}^{(0,1)})$.    By Lemma 3.3 in Crown~\cite{cr}, $$\dim HC_{r}(\Delta(E_n^{(1,e_1)})) = \binom{n-k-1}{r}.$$
Since the set of homology representatives of $\Delta(T_{n-k+1}^{(0,1)})$ is equal to the set of homology representatives of $\Delta(E_n)$ where 1 and 2 are not in the same block, using the relabeling and the set of homology representatives for $HC(\Delta(E_n^{(0,e_1)}))$, $\alpha_r$ is injective for all $r$, $n-k-1 \geq r \geq -1$.  Therefore, since the above sequence is exact,
$$\dim HC_{n-k}(\Delta(E_n^{(0,e_2)})) = \binom{n-1}{n-k+1}$$
and
$$\dim HC_{r}(\Delta(E_n^{(0,e_2)})) = \binom{n-1}{r+1} - \binom{n-k}{r+1} - \binom{n-k-1}{r}.$$

\hspace{0.5in}We continue in this manner, successively computing the dimensions of the homology groups of $HC_{r}(\Delta(E_n^{(0,e_i)}))$.  We claim that for all $i$, the complex $\Delta(E_n^{(1,e_i)})$ is equivalent, under a relabeling of the vertices, to the complex $\Delta(E_{n-k+1}^{(0,1)})$.  In particular, we relabel vertex $i$ in $E_{n-k+1}$ with $e_i$, vertex $i+1$ with $k+i$, ..., and vertex $n-k+1$ with $n$.  The sets $\Delta_r(E_n^{(1,e_i)})$ consist of all partitions $[B_1,...,B_{r+2}]$ such that $1 \in B_1$, $e_i \in B_j$ for some $j$, and none of the edges $e_1, e_2, ..., e_{i-1}$ are contained in a block of the partition.  Notice that since $k > \lceil{\frac{n}{2}}\rceil$, the condition that vertices $i-1$ and $i$ are not in the same block of a partition in $\Delta_r(E_{n-k+1}^{(0,1)})$ is equivalent to the condition that none of the edges $e_1, e_2, ... e_{i-1}$ are contained in a block of a partition in $\Delta_r(E_n^{(1,e_i)})$.  Thus, for $2 \leq i \leq n-k+1$,
$$\dim HC_{r}(\Delta(E_n^{(1,e_i)})) = \binom{n-k-1}{r}.$$
By a similar argument to the one used above, one can then show that the maps $\alpha_r$ are injective, and by exactness of each of the long exact sequences, we have the desired result. 

\qed

We are now able to compute the dimensions of the homology groups of the cyclic coloring of a diagonal hypergraph.

\begin{theorem}Let $H$ be a diagonal hypergraph with $k > \lceil{\frac{n}{2}}\rceil$.  For $n-k-1 \geq r \geq 0$,
$$\dim HC_{r}(\Delta(H)) = (n-k)\binom{n-k-1}{r} + \binom{n-k}{r+1},$$
and
$$\dim HC_{-1}(\Delta(H)) = 1.$$
\end{theorem}

\PROOF

Consider the long exact sequence
\begin{center}
\begin{tabular}{cccccccc}
0 & $\rightarrow$ & $HC_{n-k}(\Delta(E_n))$ & $\rightarrow$ & $HC_{n-k}(\Delta(H)^C)$ & $\stackrel{\phi_{n-k}}{\rightarrow}$ & \\
$HC_{n-k-1}(\Delta(H))$ & $\stackrel{\alpha_{n-k-1}}{\rightarrow}$ & $HC_{n-k-1}(\Delta(E_n))$ & $\stackrel{\beta_{n-k-1}}{\rightarrow}$ & $HC_{n-k-1}(\Delta(H)^C)$ & $\stackrel{\phi_{n-k-1}}{\rightarrow}$ & \\
$HC_{n-k-2}(\Delta(H))$ & $\stackrel{\alpha_{n-k-2}}{\rightarrow}$ & $HC_{n-k-2}(\Delta(E_n))$ & $\stackrel{\beta_{n-k-2}}{\rightarrow}$ & $HC_{n-k-2}(\Delta(H)^C)$ & $\stackrel{\phi_{n-k-2}}{\rightarrow}$ & $\hdots$ \\
\end{tabular}
\end{center}
From the proof of the last theorem, we know that the set of homology representatives of $\Delta(HC_{r}(\Delta(H)^C))$ is a subset of the set of homology representatives of $HC_{r}(\Delta(E_n))$.  Thus, the dimension of the image $\beta_{r}$ equals the dimension of $HC_{r}(\Delta(H)^C)$.  This implies that the dimension of the kernel of $\phi_r$ is zero for all $r$.  By exactness, the theorem then follows.
\qed

\bibliography{biblio}

\end{document}